\documentclass{article}
\usepackage{graphicx} 
\usepackage[a4paper,top=3cm,bottom=2cm,left=2.3cm,right=2.3cm,marginparwidth=1.75cm]{geometry}
\usepackage{graphicx}
\usepackage{verbatim}
\usepackage{hyperref}       
\usepackage{url}            
\usepackage{booktabs}       
\usepackage{amsfonts}       
\usepackage{nicefrac}       
\usepackage{microtype}      
\usepackage{amsmath}
\usepackage{multirow}
\usepackage{amsthm}
\usepackage{amssymb}

\usepackage{subfigure}
\usepackage{times}
\usepackage{color}

\newcommand{\bx}{\boldsymbol{x}}

\newcommand{\by}{\boldsymbol{y}}

\newcommand{\bv}{\boldsymbol{v}}

\usepackage{algorithm}
\usepackage{algpseudocode}
\title{Tackling the Curse of Dimensionality in Fractional and Tempered Fractional PDEs with Physics-Informed Neural Networks}
\author{Zheyuan Hu\thanks{Department of Computer Science, National University of Singapore, Singapore, 119077 (\href{mailto:e0792494@u.nus.edu}{e0792494@u.nus.edu},\href{mailto:kenji@nus.edu.sg}{kenji@nus.edu.sg})} \and Kenji Kawaguchi\footnotemark[1] \and Zhongqiang Zhang\thanks{Department of Mathematical Sciences, Worcester Polytechnic Institute, Worcester, MA 01609 USA (\href{mailto:zzhang7@wpi.edu}{zzhang7@wpi.edu})}\and  George Em Karniadakis\thanks{Division of Applied Mathematics, Brown University, Providence, RI 02912, USA (\href{mailto:george\_karniadakis@brown.edu}{george\_karniadakis@brown.edu})} \ \thanks{Advanced Computing, Mathematics and Data Division, Pacific Northwest National Laboratory, Richland, WA, United States}}
\date{}

\begin{document}

\maketitle

\begin{abstract}
Fractional and tempered fractional partial differential equations (PDEs) are effective models of long-range interactions, anomalous diffusion, and non-local effects. Traditional numerical methods for these problems are mesh-based, thus struggling with the curse of dimensionality (CoD).  Physics-informed neural networks (PINNs) offer a promising solution due to their universal approximation, generalization ability, and mesh-free training.
In principle, Monte Carlo fractional PINN (MC-fPINN) estimates fractional derivatives using Monte Carlo methods and thus could lift CoD. However, this may cause significant variance and errors, hence affecting convergence; in addition,  MC-fPINN is sensitive to hyperparameters. 
In general, numerical methods and specifically PINNs for tempered fractional PDEs are under-developed.
Herein, we extend MC-fPINN to tempered fractional PDEs to address these issues, resulting in the Monte Carlo tempered fractional PINN (MC-tfPINN). To reduce possible high variance and errors from Monte Carlo sampling, we replace the one-dimensional (1D) Monte Carlo with 1D Gaussian quadrature, applicable to both MC-fPINN and MC-tfPINN.
We validate our methods on various forward and inverse problems of fractional and tempered fractional PDEs, scaling up to 100,000 dimensions. 
Our improved MC-fPINN/MC-tfPINN using quadrature consistently outperforms the original versions in accuracy and convergence speed in very high dimensions.
\end{abstract}

\section{Introduction}
Fractional partial differential equations (PDEs) extend the concept of differentiation to non-integer orders, allowing more flexibility in modeling complex phenomena. They have found applications in various fields, such as physics, biology, finance, and engineering.
Fractional PDEs usually involve time-fractional derivatives or fractional Laplacians. These non-integer order operators capture non-stadard behavior such as long-range interactions, anomalous diffusion, and scale invariance, which are not adequately described by integer-order derivatives.
Furthermore, tempered fractional PDEs introduce a tempering parameter that relaxes the long-range interactions of the fractional operator. This tempering parameter allows a transition between standard integer-order PDEs and fully non-local fractional PDEs.
Tempered fractional PDEs are particularly useful for modeling phenomena of a balance between local interactions and non-local effects. By tuning the tempering parameter, one can adjust the influence of non-locality in the system, making tempered fractional PDEs versatile tools for capturing a wide range of behaviors in various applications.

Given their importance in science and engineering, numerous methods have been proposed to obtain their numerical solutions, see e.g. in \cite{d2020numerical,Karniadakis-b19} and many others. 
In addition to traditional mesh-based methods,  Monte Carlo methods are also developed for fractional operators, e.g., \cite{leonenko2022monte,sheng2023efficient}. 
The physics-informed neural network (PINN) \cite{raissi2019physics} has emerged as a mesh-free approach, benefiting primarily from the powerful generalization capabilities of neural networks \cite{kawaguchi2017generalization}, their universal approximation properties \cite{barron1993universal}, and robust optimization using stochastic gradients \cite{kingma2014adam}. 
Recently, PINN has been validated as highly effective in solving high-dimensional PDEs \cite{hu2023tackling,hu2023hutchinson,hu2024score}, lifting the curse of dimensionality. Regarding fractional PDEs, fractional PINN (fPINN) \cite{pang2019fpinns} and Monte Carlo fPINN (MC-fPINN) \cite{guo2022monte} represent seminal works for solving fractional PDEs with PINNs. Despite demonstrating high accuracy and scalability across various low-dimensional and high-dimensional fractional PDE forward/inverse problems, their application to particularly high dimensions, such as problems in thousands of dimensions, remains limited. Moreover, the PINN approach is lacking for the equally significant and more practical tempered fractional PDEs.

Furthermore, fPINN \cite{pang2019fpinns} suffers from the curse of dimensionality as it uses a grid-based method to discretize the fractional derivatives. 
Although MC-fPINN \cite{guo2022monte} employs Monte Carlo to overcome this issue, the computation of the fractional derivative's singularity-integrated integral using Monte Carlo may lead to significant variance and errors, hindering PINN convergence.
Moreover, MC-fPINN relies on many hyperparameters, such as setting a radius to divide the integral of the fractional derivative into internal and external parts, and it further requires a truncation threshold for the former integral to avoid numerical overflow.
Therefore, the accuracy of MC-fPINN is highly sensitive to these hyperparameters, thereby limiting its practical applicability. Additionally, while recent efforts have scaled up PINN to thousands of dimensions \cite{hu2023tackling,hu2023hutchinson}, these endeavors primarily address integer-order PDEs and do not apply to the fractional-order PDEs considered in this work.

To address these issues, this paper makes two core contributions: (1) We extend MC-fPINN \cite{guo2022monte} to tempered fractional PDEs, resulting in Monte Carlo tempered fractional PINN, denoted as MC-tfPINN. (2) Furthermore, to overcome the high variance and error associated with Monte Carlo in MC-fPINN, as well as its dependency on many hyperparameters, we propose replacing the 1D Monte Carlo computation for singularity-integrated parts with a 1D Gaussian quadrature. This quadrature method is also applied to MC-tfPINN, resulting in improved versions.

Specifically, MC-fPINN primarily relies on Monte Carlo sampling from a Beta distribution. For tempered fractional PDEs, we use Monte Carlo sampling from a Gamma distribution. Additionally, due to the practicality introduced by the tempering effect, MC-tfPINN does not need to split the integral of the tempered fractional derivative into two parts for calculation, making its algorithm more concise.

For the second part of improving MC-fPINN and MC-tfPINN, we observe that simulating a $d$-dimensional fractional PDE operator requires computing a $d$-dimensional integral, which can be reformulated as a product integral of a 1D radial singular integral and an integral over a $d$-dimensional unit sphere. 
For the former, the singularity of the integral introduces significant variance and error to Monte Carlo, leading us to propose using Gaussian quadrature to compute this 1D radial singular integral. 
Compared to Monte Carlo sampling, integration using the Gauss quadrature provides theoretical guarantees of accuracy when computing singular integrals, making it more precise. Additionally, with Monte Carlo, resampling at the beginning of each optimization epoch is required, which requires additional time, whereas quadrature points remain fixed and do not require resampling. 
Moreover, as quadrature points are fixed, we do not need a truncation threshold to avoid numerical overflow.
In the Monte Carlo method, such a threshold is necessary due to the randomness of the sampling to prevent sampling extreme values.
Furthermore, regarding MC-fPINN's need to set a radius to divide the integral of the fractional derivative into internal and external parts, we choose this radius as the diameter of the PDE support domain, thereby reducing the dependence on hyperparameters and simplifying the computation in the second part, thus further accelerating PINN.

We validate the effectiveness of our proposed methods on several fractional and tempered fractional PDEs, scaling up to 100,000 dimensions, far exceeding previous works. 
Firstly, for the fractional Poisson equation, we test a complex anisotropic solution, confirming the higher accuracy and faster training speed of the improved MC-fPINN over the regular MC-fPINN. 
%
%
Next, we compare MC-tfPINN and the quadrature-based improved MC-tfPINN on tempered fractional Poisson and diffusion equations, once again confirming better efficiency and higher accuracy of our methods. 
Additionally, we vary the fractional order and tempering coefficient in the tempered fractional PDE operator, demonstrating the stability of our models across various dimensions. 
Finally, we test inverse problems and show that our proposed models are versatile and applicable.

In summary,  the novelty and contribution of the work lies in the following aspects: 
1) extending MC-fPINN to tempered fractional PDEs; 
and 2)introducing improved versions based on quadrature over Monte Carlo for fPINN, rendering simulations of high-dimensional fractional and tempered fractional PDEs more accurate and efficient. 

\section{Related Work}
For fractional PDE solvers, Pang et al. \cite{pang2019fpinns} introduce a fractional Physics-Informed Neural Network (fPINN), utilizing neural networks as proxies for solutions and employing a discretized approach to fractional derivatives for training. However, the discretization in fPINN is grid-based and thus suffers from the curse of dimensionality.
Guo et al. \cite{guo2022monte} develop a Monte Carlo fractional PINN (MC-fPINN), which approximates Caputo-type time-fractional derivatives and fractional Laplacian within a hyper-singular integral framework, leveraging Monte Carlo methods based on Beta distributions.
Ma et al. \cite{ma2023bi} further propose 	
bi-orthogonal fPINN to deal with stochastic fractional PDEs.
Chen et al. \cite{chen2021solving} use PINNs to address both forward and inverse problems associated with Fokker-Planck-Lévy equations and use the discretization scheme in \cite{gao2016fokker} to estimate fractional derivatives.
Firoozsalari et al. \cite{firoozsalari2023deepfdenet} employ Gaussian quadrature for calculating numerical integrals associated with fractional derivatives, concentrating exclusively on time derivatives and neglecting to address scenarios involving high-dimensional fractional Laplacians.
Leonenko and Podlubny \cite{leonenko2022monte} introduce a Monte Carlo method to estimate the Grünwald–Letnikov fractional derivative.
Based on the approach by Leonenko and Podlubny \cite{leonenko2022monte}, Wang and Karniadakis \cite{wang2024gmc} propose a general quasi Monte Carlo PINN for solving fractional PDEs on irregular domains.
Sheng et al. \cite{sheng2023efficient} propose a novel Monte Carlo method for fractional Laplacian in high dimensions.

For tempered fractional PDEs, a pedagogical guide for tempered fractional calculus is provided in \cite{SABZIKAR201514}. The pedagogical guide for the tempered factional diffusion process is also presented in \cite{BAEUMER20102438}, where the tempered stable distribution is discussed for the tempered Lévy process.
Duo et al. \cite{duo2019numerical} propose an accurate finite difference method for tempered fractional Laplacian and research the tempered effect.
Sun et al. \cite{sun2021algorithm} provide algorithm implementation and analysis for 2D tempered fractional Laplacian.
Deng et al. \cite{deng2018boundary} solve boundary problems of fractional and tempered fractional PDEs.
Li and Deng \cite{li2016high} propose several high-order numerical schemes for the tempered fractional diffusion equations.
Hanert and Piret \cite{hanert2014chebyshev} propose a Chebyshev pseudospectral method for the tempered fractional diffusion equations.
There also exists a spectral method for tempered fractional cases \cite{zhao2016spectral}.

\section{Preliminary}

We consider high-dimensional tempered fractional PDEs on a bounded spatial domain $\Omega \subset \mathbb{R}^d$:
\begin{align} \label{eq:pde-tempered-fractional}
\frac{\partial^{\gamma, \lambda_t} u(\bx, t)}{\partial t^\gamma} + c \left(-\Delta_{\lambda_{\bx}}\right)^{\alpha/2}u(\bx, t) + \bv \cdot \nabla_{\bx}u(\bx, t) &= f(\bx, t), \quad (\bx, t) \in \Omega \times (0, T].\\
u(\bx, t=0)&= g(\bx), \quad \bx \in \Omega.\\
u(\bx, t) &= 0, \quad (\bx, t) \in \Omega^c \times (0, T].
\end{align}
Here, $f(\bx, t)$ is the forcing term, $c \in \mathbb{R}^+$ is the diffusion coefficient, $\lambda_t\geq 0$ and $\lambda_{\bx}\geq 0$ are the tempering coefficients, and $\bv \in \mathbb{R}^d$ is the flow velocity.
The time-tempered fractional derivative is defined as \cite{fahad2021tempered}:
\begin{align}\label{eq:def_tempered_fractional_time_derivative}
\frac{\partial^{\gamma,\lambda} f(t)}{\partial t^\gamma} := \frac{1}{\Gamma(1 - \gamma)} \int_0^t \exp\left(-\lambda(t - \tau)\right)(t - \tau)^{-\gamma}\frac{\partial f(\tau)}{\partial \tau} d\tau =e^{-\lambda t} \frac{\partial^{\gamma} \left(e^{\lambda t} f(t)\right)}{\partial t^\gamma}, \quad 0 < \gamma < 1.
\end{align}
When $\lambda=0$, the tempered fractional derivative becomes the Caputo time-fractional derivative. 
The tempered fractional Laplacian is defined \cite{deng2018boundary}:
\begin{align}\label{eq:fractional-laplacian-tempered}
(-\Delta_\lambda)^{\alpha / 2} u :
= C_{d,\alpha,\lambda} \text{P.V.} \int_{\mathbb{R}^{d}} \frac{u(\bx) - u(\by)}{\exp\left(\lambda\Vert\bx - \by\Vert\right)\Vert \bx - \by \Vert^{d + \alpha}}d\by, \quad \alpha \in (0, 2),
\end{align}
where P.V. denotes the principle value of the integral and $C_{d,\alpha,\lambda} = \frac{\Gamma\left(\frac{d}{2}\right)}{2\pi^{\frac{d}{2}}\left|\Gamma\left(-{\alpha}\right)\right|}$.
Note that when $\lambda=0$, the tempered fractional Laplacian becomes the fractional Laplacian e hyper-singular integral \eqref{eq:fractional-laplacian}. 
When $\lambda_t=\lambda_{\bx}=0$, the equation \eqref{eq:pde-tempered-fractional} becomes the fractional PDEs considered in \cite{guo2022monte}.
For simplicity, we denote 
\begin{align}
\mathcal{L}^\phi\left[u(\bx, t)\right] :=
\frac{\partial^{\gamma, \lambda_t} u(\bx, t)}{\partial t^\gamma} + c \left(-\Delta_{\lambda_{\bx}}\right)^{\alpha/2}u(\bx, t) + \bv \cdot \nabla_{\bx}u(\bx, t),
\end{align}
and $\phi = \{\gamma, c, \alpha, \bv\}$ is a collection of coefficients and parameters in the problem above. 
\subsection{Tempered fractional PDEs: forward and inverse problems}
We consider forward and inverse tempered fractional PDE problems.
In the forward problem, we aim to solve the PDE solution $u(\bx, t)$
given all coefficients/parameters $\phi$), the forcing term $f(\bx, t)$ and initial and boundary conditions.
In the inverse problem, we have sparse observation data over the exact solution $u$ while we are given initial and boundary conditions. The goal is to find the unknown coefficients/parameters in the PDE operator, e.g., $\lambda, \alpha, \gamma, \boldsymbol{v}$, etc.


\subsection{Physics-Informed Neural Networks (PINNs)}
We use PINNs \cite{raissi2019physics} to solve both forward and inverse problems of the fractional and tempered fractional PDEs defined in the previous subsections. 

%
PINN is a neural network-based PDE solver that parameterizes the PDE surrogate solution using the model $u_\theta(\bx, t)$ with parameters $\theta$. PINN optimizes the neural network model 
$u_\theta(\bx, t)$ to satisfy the PDE operator, initial, and boundary conditions. As long as PINN approximates these three conditions well enough, it can approximate the PDE solution accurately. See PINN's theoretical guarantees and foundations in \cite{de2022error,hu2021extended,Luo2020TwoLayerNN,mishra2020estimates}.

\textbf{Boundary condition}. We hard constraint the vanishing boundary condition using a function $\rho(\bx)$, such that $\rho(\bx) > 0$ for $\bx \in \Omega$ and $\rho(\bx) = 0$ for $\bx \in \Omega^c$. 
i.e., the surrogate model will be
\begin{align}
u_\theta(\bx, t) = \rho (\bx) \cdot u_\theta^{\text{NN}}(\bx, t),
\end{align}
where $u_\theta^{\text{NN}}$ is a full-connected network and $u_\theta$ will be the final surrogate model.
Here, hard constraints are essential for stable training and successful optimization of PINNs.

\textbf{Initial condition}. Given the initial points $\{\bx_{n,\text{initial}}\}_{n=1}^{N_{\text{initial}}} \subset \Omega$, we minimize the $L_2$ distance between PINN prediction and exact initial condition:
\begin{align}
\mathcal{L}_{\text{initial}}(\theta) = \frac{1}{N_{\text{initial}}}\sum_{n=1}^{N_{\text{initial}}} {|u_{\theta}(\bx_{n,\text{initial}}, 0)-g(\bx)|}^2.
\end{align}

\textbf{PDE operator}. Given the residual/collocation points $\{\bx_{n,\text{residual}}, t_{n,\text{residual}}\}_{n=1}^{N_{\text{residual}}} \subset \Omega \times (0, T]$, we minimize the PINN residual:
\begin{equation}
\mathcal{L}_{\text{residual}}(\theta, \phi)=\frac{1}{N_{\text{residual}}}\sum_{n=1}^{N_{\text{residual}}} |\mathcal{L}^\phi u_{\theta}(\bx_{n,\text{residual}}, t_{n,\text{residual}})-f(\bx_{n,\text{residual}}, t_{n,\text{residual}})|^2.
\end{equation}
In forward problems, we have full knowledge of PDE coefficients $\phi$, and thus, we only optimize $\theta$ while in inverse problems, we optimize both PINN model parameters $\theta$ and PDE unknown coefficients $\phi$.

\textbf{Data loss}. Inverse problems allow access to sparse observation over the exact solution $\{\bx_{n, \text{data}}, t_{n, \text{data}}, u^*_{n}\}_{n=1}^{N_{\text{data}}}$, which further leads to the following data loss for the PINN surrogate model to capture these observation data to recover the unknown PDE coefficients $\phi$:
\begin{align}
\mathcal{L}_{\text{data}}(\theta) = \frac{1}{N_{\text{data}}}\sum_{n=1}^{N_{\text{data}}} {|u_{\theta}(\bx_{n, \text{data}}, t_{n, \text{data}})-u^*_n|}^2.
\end{align}

In addition to these loss functions, we usually weighted sum them to form the final loss. The weights are denoted $\lambda_{\text{initial}}$ for the initial loss, $\lambda_{\text{residual}}$ for the residual loss, and $\lambda_{\text{data}}$ for the data loss, respectively.
Given the different backgrounds of forward and inverse problems, we summarize their respective loss functions for PINNs.
\begin{itemize}
\item In the forward problem, we optimize the initial and residual losses:
\begin{align}
\min_\theta \lambda_{\text{initial}}\mathcal{L}_{\text{initial}}(\theta) + \lambda_{\text{residual}}\mathcal{L}_{\text{residual}}(\theta).
\end{align}
\item In the inverse problem, we optimize all losses:
\begin{align}
\min_{\theta,\phi}\lambda_{\text{initial}}\mathcal{L}_{\text{initial}}(\theta) + \lambda_{\text{residual}}\mathcal{L}_{\text{residual}}(\theta,\phi) +\lambda_{\text{data}} \mathcal{L}_{\text{data}}(\theta).
\end{align}
\end{itemize}

\subsection{Monte Carlo fPINN (MC-fPINN)}
Although PINN, combined with modern automatic differentiation libraries, can utilize the above loss function to solve integer-order PDEs, automatic differentiation does not cover fractional derivatives.
In Guo et al. \cite{guo2022monte}, Monte Carlo fPINN is proposed to simulate the fractional Laplacian and derivatives to enable PINN's application to fractional PDEs.

\subsubsection{Fractional Laplacian}
In Monte Carlo-based approaches, e.g., Monte-Carlo fPINN \cite{guo2022monte}, the following fractional Laplacian is defined as the
hyper-singular integral \cite{LISCHKE2020109009}:
\begin{align}\label{eq:fractional-laplacian}
(-\Delta)^{\alpha / 2} u :
= C_{d,\alpha} \text{P.V.} \int_{\mathbb{R}^{d}} \frac{u(\bx) - u(\by)}{\Vert \bx - \by \Vert^{d + \alpha}}d\by, \quad \alpha \in (0, 2),
\end{align}
where P.V. denotes the principle value of the integral and $C_{d,\alpha} = \frac{2^\alpha\Gamma\left(\frac{d + \alpha}{2}\right)}{\pi^{\frac{d}{2}}\left|\Gamma\left(-\frac{\alpha}{2}\right)\right|}$.
According to \cite{guo2022monte}, the integral is divided into two parts for Monte Carlo simulation, 
\begin{align}
(-\Delta)^{\alpha / 2} u
= C_{d,\alpha} \left(\int_{\by \in B_{r_0}(\bx)} \frac{u(\bx) - u(\by)}{\Vert \bx - \by \Vert^{d + \alpha}}d\by + \int_{\by \notin B_{r_0}(\bx)} \frac{u(\bx) - u(\by)}{\Vert \bx - \by \Vert^{d + \alpha}}d\by\right),
\end{align}
where $B_{r_0}(\bx) = \{\by \ |\ \Vert \by - \bx \Vert \leq r_0\}$.  \cite{guo2022monte} shows that Monte Carlo can simulate the two parts:
\begin{align}
\int_{\by \in B_{r_0}(\bx)} \frac{u(\bx) - u(\by)}{\Vert \bx - \by \Vert^{d + \alpha}}d\by = \frac{|\mathcal{S}^{d-1}|r_0^{2-\alpha}}{2(2-\alpha)}\mathbb{E}_{\boldsymbol{\xi},r \sim f_1(r)}\left[\frac{2u(\bx) - u(\bx - \boldsymbol{\xi} r) - u(\bx + \boldsymbol{\xi} r)}{r^{2}}\right],
\end{align}
where we sample $r$ based on $r / r_0 \sim \operatorname{Beta}(2-\alpha,1)$ and $\boldsymbol{\xi}$ uniformly from the unit sphere. To prevent the numerical instability due to $r \approx 0$, \cite{guo2022monte} adopts a truncation to ensure the sampled $r$ is larger than some fixed positive $\epsilon$. 
The second integral defined outside of the ball $B_{r_0}$ is rewritten as 
\begin{align}
\int_{\by \notin B_{r_0}(\bx)} \frac{u(\bx) - u(\by)}{\Vert \bx - \by \Vert^{d + \alpha}}d\by = \frac{|\mathcal{S}^{d-1}|r_0^{-\alpha}}{2\alpha}\mathbb{E}_{\boldsymbol{\xi},r \sim f_2(r)}\left[2u(\bx) - u(\bx - \boldsymbol{\xi} r) - u(\bx + \boldsymbol{\xi} r)\right],
\end{align}
where we sample $r$ based on $r_0 / r \sim \operatorname{Beta}(\alpha,1)$ and $\boldsymbol{\xi}$ uniformly from the unit sphere.

\subsubsection{Time-fractional Derivative}
Monte Carlo on the Beta distribution can also evaluate the time-fractional derivative,  see Guo et al. \cite{guo2022monte} or  \eqref{eq:time-fractional-derivative-equiv} for the detailed derivation:
\begin{align}\label{eq:time-fractional-derivative-as-expecation}
\frac{\partial^\gamma f(t)}{\partial t^\gamma} = \frac{\gamma}{1 - \gamma}t^{1-\gamma}\mathbb{E}_{\tau \sim \text{Beta}(1-\gamma, 1)}\left[\frac{f(t) - f(t - t\tau)}{t\tau}\right]+ \frac{f(t) - f(0)}{t^\gamma}.
\end{align}
Hence, we can simulate the differential operators in the fractional PDEs using Monte Carlo simulation. Then, we can plug these PDE operator simulators into the PINN loss to enable PINN to solve them.

\section{Proposed Method}
This section presents the main proposed methods. Specifically, we first extend MC-fPINN \cite{guo2022monte} to tempered fractional cases in subsection \ref{sec:mc_tfpinn} using a pure Monte Carlo-based approach, enabling the simulation of tempered fractional Laplacian to enable PINN's application to tempered fractional PDEs.
Then, in subsection \ref{sec:improved_mc_fpinn} we introduce our improved MC-fPINN that utilizes a more accurate and efficient quadrature to estimate the 1D integral about $r$ in the fractional derivatives.
Finally, subsection \ref{sec:improved_mc_tfpinn} extends the improving techniques to tempered fractional cases.
%

\subsection{Monte Carlo Tempered fPINN (MC-tfPINN)}\label{sec:mc_tfpinn}
We first extend MC-fPINN \cite{guo2022monte} to tempered fractional cases using a Monte Carlo-based approach.
The tempered fractional Laplacian can be similarly defined \cite{deng2018boundary}:
\begin{align}
(-\Delta_\lambda)^{\alpha / 2} u :
= C_{d,\alpha,\lambda} \text{P.V.} \int_{\mathbb{R}^{d}} \frac{u(\bx) - u(\by)}{\exp\left(\lambda\Vert\bx - \by\Vert\right)\Vert \bx - \by \Vert^{d + \alpha}}d\by, \quad \alpha \in (0, 2),
\end{align}
where P.V. denotes the principle value of the integral and $C_{d,\alpha,\lambda} = \frac{\Gamma\left(\frac{d}{2}\right)}{2\pi^{\frac{d}{2}}\left|\Gamma\left(-{\alpha}\right)\right|}$.
First, we represent the fractional Laplacian as an expectation:
\begin{align}\label{eq:fractional-laplacian-as-expectation}
\int_{\mathbb{R}^{d}} \frac{u(\bx) - u(\by)}{\exp\left(\lambda\Vert\bx - \by\Vert\right)\Vert \bx - \by \Vert^{d + \alpha}}d\by &= \int_{\mathbb{R}^{d}} \frac{u(\bx) - u(\bx - \by)}{\exp\left(\lambda\Vert\by\Vert\right)\Vert \by \Vert^{d + \alpha}}d\by \notag \\
&= \frac{1}{2} \int_{\mathbb{R}^{d}} \frac{2u(\bx) - u(\bx - \by) - u(\bx + \by)}{\exp\left(\lambda\Vert\by\Vert\right)\Vert \by \Vert^{d + \alpha}}d\by \notag  \\
&= \frac{1}{2} \int_{\mathcal{S}^{d-1}}\int_{0}^\infty \frac{2u(\bx) - u(\bx - \boldsymbol{\xi} r) - u(\bx + \boldsymbol{\xi} r)}{\exp\left(\lambda r\right)r^{1 + \alpha}}drd\boldsymbol{\xi}\notag \\
&= \frac{1}{2} \int_{\mathcal{S}^{d-1}}\int_{0}^\infty \frac{2u(\bx) - u(\bx - \boldsymbol{\xi} r) - u(\bx + \boldsymbol{\xi} r)}{r^{2}}\exp\left(-\lambda r\right)r^{1-\alpha}drd\boldsymbol{\xi} \notag \\
&= \frac{1}{2}|\mathcal{S}^{d-1}|\frac{\Gamma(2 - \alpha)}{\lambda^{2-\alpha}} \mathbb{E}_{r, \boldsymbol{\xi}} \left[\frac{2u(\bx) - u(\bx - \boldsymbol{\xi} r) - u(\bx + \boldsymbol{\xi} r)}{r^{2}}\right],
\end{align}
where $r$ obeys the Gamma distribution $\operatorname{Gamma}(2 - \alpha, \lambda)$ and $\boldsymbol{\xi}$ is uniformly distributed on the $d$-dimensional unit sphere.

Similar to the Monte Carlo fPINN case \cite{guo2022monte}, 
we use Monte Carlo methods to simulate the expectation over $r$ and 
$\boldsymbol{\xi}$. Also, we pick the sampled $r\geq \epsilon>0$  to prevent the numerical instability due to $r \approx 0$. 
When $\lambda >0$, thanks to the tempered term $\exp(-\lambda r)$, we don't need to decompose the sampling into two intervals like Monte Carlo fPINN \cite{guo2022monte}, reducing the hyperparameter number, simplifying the algorithm and stabilizing the training. Besides, as $r \rightarrow \infty$, the sampled Gamma distribution exhibits a short tail, reducing the variance. 

We have introduced the Monte Carlo simulation for tempered fractional Laplacian. Meanwhile, the approximation of the tempered time factional derivative can be reduced to that of the conventional fractional case; see equation (\ref{eq:def_tempered_fractional_time_derivative}). Hence, the methodology will be the same as MC-fPINN for the tempered time factional derivative.

\subsection{Improved MC-fPINN}\label{sec:improved_mc_fpinn}
In this subsection, we propose a new class of improved MC-fPINN to accelerate training and achieve lower error based on quadrature.
%
MC-fPINN suffers from two issues. Firstly, MC-fPINN relies heavily on numerous hyperparameters, particularly $r_0$ to divide the integral within the fractional Laplacian into two parts, and $\epsilon$ to mitigate numerical instability, rendering the model more sensitive to hyperparameter tuning. For instance, with $\epsilon$, setting it too small leads to numerical instability and non-convergence, while setting it too large results in bias and error in estimating the fractional Laplacian.

To this end, we propose an improved MC-fPINN to reduce its hyperparameter sensitivity. We choose $r_0$ as the diameter of the PDE's support domain $\Omega$ by default.
Then, we can use the Gauss-Jacobi quadrature to compute the integral concerning $r$ from 0 to $r_0$ with a singularity at $r=0$:
\begin{equation}
\int_{-1}^1 (1-x)^\alpha (1+x)^\beta f(x)dx \approx \sum_{i=0}^{N_{\text{quad}}}w_i f(x_i).
\end{equation}
We rescale the integral interval, choose $\alpha,\beta$ accordingly, and rescale to perform the Gauss-Jacobi quadrature. Then, we compute the 1D integral on $r$ via quadrature and the $(d-1)$-dimensional integral on the unit sphere via Monte Carlo.
\begin{align}
\int_{\by \in B_{r_0}(\bx)} \frac{u(\bx) - u(\by)}{\Vert \bx - \by \Vert^{d + \alpha}}d\by &=\frac{1}{2} \int_{\mathcal{S}^{d-1}}\int_{0}^{r_0} \frac{2u(\bx) - u(\bx - \boldsymbol{\xi} r) - u(\bx + \boldsymbol{\xi} r)}{r^{2}}r^{1-\alpha}drd\boldsymbol{\xi}\\
&\approx \frac{|\mathcal{S}^{d-1}|}{2} \sum_{i=1}^Nw_i\frac{2u(\bx) - u(\bx - \boldsymbol{\xi}_i r_i) - u(\bx + \boldsymbol{\xi}_i r_i)}{r_i^{2}},
\end{align}
where $(w_i, r_i)_{i=1}^N$ are Gauss-Jacobi quadrature weight and points, while $\boldsymbol{\xi}_i$ are sampled uniformly from the unit sphere.
For the second part of the fractional Laplacian, note that $u(\bx \pm \boldsymbol{\xi} r)$ will be 0 since $r \geq r_0$ and $r_0$ is the diameter of the support $\Omega$:
\begin{align}
\int_{\by \notin B_{r_0}(\bx)} \frac{u(\bx) - u(\by)}{\Vert \bx - \by \Vert^{d + \alpha}}d\by = \frac{|\mathcal{S}^{d-1}|r_0^{-\alpha}}{\alpha}u(\bx),
\end{align}
The Gauss-Jacobi quadrature for the integral about $r$ has multiple advantages.
\begin{itemize}
\item We do not require two hyperparameters, $r_0$ and $\epsilon$. 
For the integral about $r$, we replace the Monte Carlo sampling for $r$ with a more precise quadrature. This is far more accurate than Monte Carlo, especially when dealing with integrals that have singularities. 
\item  No resampling over $r$ is performed as quadrature points are fixed.
In the regular MC-fPINN, each iteration requires sampling new $r$ values from the Beta distribution. Thus, our improved method saves time when sampling.
\item The second part of the integral in our method is simplified and requires less neural network inference, effectively reducing the overall computational load by about 50\%.
\end{itemize}

A similar quadrature can evaluate the time-fractional derivative. We use a Gauss-Jacobi quadrature for the integral:
\begin{align}
\frac{\partial^\gamma f(t)}{\partial t^\gamma} &= \frac{1}{\Gamma(1 - \gamma)} \int_0^t (t - \tau)^{-\gamma}\frac{\partial f(\tau)}{\partial \tau} d\tau \notag \\
&=\gamma \int^t_0 \tau^{-\gamma} \frac{f(t) - f(t - \tau)}{\tau}d\tau + \frac{f(t) - f(0)}{t^\gamma}\notag \\
&=\gamma \int^1_0 (t\tau)^{-\gamma} \frac{f(t) - f(t - t\tau)}{t\tau}dt\tau + \frac{f(t) - f(0)}{t^\gamma}\notag\\
&=\gamma t^{1-\gamma}\int^1_0 (\tau)^{-\gamma} \frac{f(t) - f(t - t\tau)}{t\tau}d\tau + \frac{f(t) - f(0)}{t^\gamma} \label{eq:time-fractional-derivative-equiv} \\
&\approx \gamma t^{1-\gamma}\left(\sum_{i=1}^N w_i \frac{f(t) - f(t - t\tau_i)}{t\tau_i}\right) + \frac{f(t) - f(0)}{t^\gamma}, \label{eq:time-fractional-derivative-quad}
\end{align}
where $(w_i, \tau_i)_{i=1}^N$ are Gauss-Jacobi quadrature weight and points.

\subsection{Improved tfPINN}\label{sec:improved_mc_tfpinn}

For tfPINN, we estimate the integral on $r$ with generalized Gauss–Laguerre quadrature, which is tailored for the following type of integral with a singularity at $x=0$ for $\alpha > -1$ and $\lambda > 0$:
\begin{equation}
\int_0^\infty x^\alpha e^{-\lambda x} f(x)dx \approx \sum_{i=0}^{N_{\text{quad}}}w_i f(x_i).
\end{equation}
We compute the 1D integral on $r$ via quadrature and the integral on the unit sphere via Monte Carlo. By the equivalent form \eqref{eq:fractional-laplacian-as-expectation}, we write
\begin{align}
\int_{\mathbb{R}^{d}} \frac{u(\bx) - u(\by)}{\exp\left(\lambda\Vert\bx - \by\Vert\right)\Vert \bx - \by \Vert^{d + \alpha}}d\by
%
&= \frac{1}{2}|\mathcal{S}^{d-1}| \mathbb{E}_{ \boldsymbol{\xi}} \left[\int_0^\infty \frac{2u(\bx) - u(\bx - \boldsymbol{\xi} r) - u(\bx + \boldsymbol{\xi} r)}{r^{2}} \exp\left(-\lambda r\right)r^{1-\alpha}dr\right]\\
&\approx \frac{1}{2}|\mathcal{S}^{d-1}|  \sum_{i=0}^{N}w_i \frac{2u(\bx) - u(\bx - \boldsymbol{\xi}_i r_i) - u(\bx + \boldsymbol{\xi}_i r_i)}{r_i^{2}}.
\end{align}
where $(w_i, r_i)_{i=1}^N$ are Gauss-Laguerre quadrature weight and points, while $\boldsymbol{\xi}_i$ are sampled uniformly from the unit sphere.

\subsection{Discussion: On Quasi-Monte Carlo}
We have employed Gauss-quadrature rules to compute the 1D integral about $r$ more accurately and efficiently. For $\boldsymbol{\xi}$, we still utilize standard Monte Carlo because $\boldsymbol{\xi} \in \mathbb{R}^d$ and quadrature will suffer from the CoD.
To further reduce the variance of Monte Carlo over $\boldsymbol{\xi} \in \mathbb{R}^d$, Quasi-Monte Carlo (QMC) is an appealing option. In Appendix \ref{exp:mc_tfpinn_qmc}, we also compare the performance of using QMC against using Monte Carlo methods. We do not observe any improvement in accuracy using QMC. 

\section{Computational Experiments}
\subsection{Fractional Poisson Equation}
Consider the fractional Poisson equation $(-\Delta)^{\alpha / 2} u = f$ in the unit ball. We will test the following exact solution:
\begin{equation}
u_{\text{exact}}(\bx) = \left(1 - \Vert \bx \Vert^2\right)^{\alpha/2}\left(c_{1, 0} + \sum_{i=1}^d c_{1,i}\bx_i\right) + \left(1 - \Vert \bx \Vert^2\right)^{1 + \alpha/2}\left(c_{2, 0} + \sum_{i=1}^d c_{2,i}\bx_i\right)
\end{equation}
where $c_{1, i}$ and $c_{2, i}$ are independently sampled from a unit Gaussian. It is a linear combination of the functions in Table \ref{tab:dyda} to form a relatively complicated and anisotropic exact solution whose forcing term can be computed analytically.
\begin{table}[htbp]
\centering
\begin{tabular}{|c|c|}\hline
$u(x)$ in the unit ball & $(-\Delta)^{\alpha/2}u(x)$ in the unit ball \\\hline
$\left(1 - \Vert \bx \Vert^2\right)^{\alpha/2}$ & $-2^\alpha\Gamma\left(\frac{\alpha}{2}+1\right)\Gamma\left(\frac{\alpha+d}{2}\right)\Gamma\left(\frac{d}{2}\right)^{-1}$ \\\hline
$\left(1 - \Vert \bx \Vert^2\right)^{1+\alpha/2}$ & $-2^\alpha\Gamma\left(\frac{\alpha}{2}+2\right)\Gamma\left(\frac{\alpha+d}{2}\right)\Gamma\left(\frac{d}{2}\right)^{-1}\left(1 - \left(1 + \frac{\alpha}{d}\right)\Vert \bx \Vert^2\right)$\\\hline
$\left(1 - \Vert \bx \Vert^2\right)^{\alpha/2}\bx_d$ & $-2^\alpha\Gamma\left(\frac{\alpha}{2}+1\right)\Gamma\left(\frac{\alpha+d}{2}+1\right)\Gamma\left(\frac{d}{2}+1\right)^{-1}\bx_d$ \\\hline
$\left(1 - \Vert \bx \Vert^2\right)^{1+\alpha/2}\bx_d$ & $-2^\alpha\Gamma\left(\frac{\alpha}{2}+2\right)\Gamma\left(\frac{\alpha+d}{2}+1\right)\Gamma\left(\frac{d}{2}+1\right)^{-1}\left(1 - \left(1 + \frac{\alpha}{d+2}\right)\Vert \bx \Vert^2\right)\bx_d$\\\hline
\end{tabular}
\caption{Fractional Laplacian of some functions, from \cite{Dyda+2012+536+555}.}
\label{tab:dyda}
\end{table}
In 
MC-fPINN \cite{guo2022monte}, a simple isotropic exact solution $\left(1 - \Vert \bx \Vert^2\right)^{1+\alpha/2}$ is tested. We extend their experiment to complicated anisotropic exact solutions and show that our improved MC-fPINN outperforms.

For a fair comparison, we use the hyperparameters from MC-fPINN \cite{guo2022monte}. 
We choose $\alpha= 1.5$ for a test. 
We use $\epsilon=1e-6$ for preventing numerical issues and truncating the sampled $r$ from Beta distribution as in MC-fPINN. 
Here $r_0$ is set to 0.25 in MC-fPINN, according to the hyperparameter study in Guo et al. \cite{guo2022monte}.

We opt for 100 residual points and 64 points for the Monte Carlo/quadrature batch size.
We randomly sampled 20K test points within the unit ball.
We optimize both models for 1 million epochs with Adam \cite{kingma2014adam}, whose initial learning rate is 1e-3 and linearly decays to zero. 
We adopt the following model structure to satisfy the zero boundary condition with hard constraint and to avoid the boundary loss \cite{lu2021physics}
$
\text{ReLU} (1 - \Vert\bx\Vert_2^2) u_\theta(\bx),
$
where $u_\theta(\bx)$ is the neural neural network model. 
The neural network model is a four-layer, fully-connected network with 128 hidden units and a hyperbolic tangent activation function.
Since we hard-constrained the boundary condition, the only loss function will be the residual loss. Hence, there is no weighting hyperparameter.

\begin{table}[htbp]
\centering
\begin{tabular}{|c|c|c|c|c|c|}
\hline
Method & Dimension & 100 & 1,000 & 10,000 & 100,000 \\ \hline
\multirow{2}{*}{MC-fPINN (Guo et al. \cite{guo2022monte})} & Speed & 261it/s & 223it/s & 77it/s & 9it/s \\ \cline{2-6} 
 & Rel. $L_2$ & 2.86E-02 & 4.99E-02 & 5.85E-02 & 8.02E-02 \\ \hline
\multirow{2}{*}{Improved MC-fPINN (Ours)} & Speed & \textbf{1092it/s} & \textbf{747it/s} & \textbf{174it/s} & \textbf{20it/s} \\ \cline{2-6} 
 & Rel. $L_2$ & \textbf{2.84E-02} & \textbf{3.68E-02} & \textbf{4.94E-02} & \textbf{7.10E-02} \\ \hline
\end{tabular}
\caption{Results for the fractional Poisson equation. Our improved MC-fPINN runs faster and performs better with slightly lower errors than the original MC-fPINN by Guo et al. \cite{guo2022monte}. Faster speed and lower error are bold.}
\label{tab:dyda_result}
\end{table}

The results are presented in Table \ref{tab:dyda_result} with different methods' speed (iteration per second) and relative $L_2$ error. Our improved MC-fPINN consistently outperforms MC-fPINN by Guo et al. \cite{guo2022monte} regarding speed and relative $L_2$ error performance. This is attributed to the Gauss quadrature for the integral about $r$ in our method, rendering it more precise compared to Monte Carlo. Additionally, there is no need to resample the integration points for $r$ at each epoch, enhancing effectiveness and efficiency compared with vanilla MC-fPINN \cite{guo2022monte}. Moreover, our improved MC-fPINN eliminates the need for additional hyperparameters like $r_0$ and $\epsilon$, making it more robust in hyperparameter selection during practical applications.
Furthermore, in our method, the second part of the integral in estimating the fractional Laplacian is simplified and it requires less neural network inference, effectively reducing the overall computational load by about 50\%, since we choose $r_0$ as the diameter of the PDE support, which is $r_0 = 2$ in this test case.

\subsection{Tempered Fractional Poisson Equation}
We consider high-dimensional tempered fractional Poisson equations $\left(-\Delta_\lambda\right)^{\alpha/2} u = f$ within the unit ball, zero boundary conditions, and various anisotropic and inseparable solutions.
The following nontrivial exact solutions are considered in related work  \cite{hu2023tackling,hu2023bias,hu2023hutchinson} to evaluate the high-dimensional PINN, specifically one with two-body interaction and another with three-body interaction:
\begin{equation}\label{eq:sol_2body}
u_{\text{two-body}}(\bx) = \text{ReLU}\left(1 - \Vert \bx \Vert_2^2\right)\left(\sum_{i=1}^{d-1}  c_i \sin(\bx_i +\cos(\bx_{i+1})+\bx_{i+1}\cos(\bx_i))\right),
\end{equation}
\begin{equation}\label{eq:sol_3body}
u_{\text{three-body}}(\bx) = \text{ReLU}\left(1 - \Vert \bx \Vert_2^2\right)\left(\sum_{i=1}^{d-2}  c_{i} \exp(\bx_i\bx_{i+1}\bx_{i+2})\right),
\end{equation}
where $c_i \sim \mathcal{N}(0, 1)$ and there are two-body interactions between each $(\bx_i, \bx_{i+1})$, three-body interactions between each pair $(\bx_i, \bx_{i+1}, \bx_{i+2})$, respectively.
The forcing terms for training residual points are computed with a Monte Carlo simulation with 1024 samples.

 We choose $\alpha = 0.5/1.5$ and $\lambda = 1$ for the tempered fractional Laplacian. We set $\epsilon = 1e-6$ for MC-tfPINN to truncate $r$ to prevent numerical instability.
The Monte Carlo/quadrature sample size for the neural network is 64 during training. 

The model is a 4-layer fully connected network with 128 hidden units activated by $\tanh$, which is trained via Adam \cite{kingma2014adam} for 10K epochs, with an initial learning rate 1e-3, which linearly decays to zero at the end of the optimization. 
We select 100 random residual points at each Adam epoch and 20K fixed testing points uniformly from the unit ball.
We adopt the following model structure to satisfy the zero boundary condition with hard constraint and to avoid the boundary loss \cite{lu2021physics}$
\text{ReLU} (1 - \Vert\bx\Vert_2^2) u_\theta(\bx),
$
where $u_\theta(\bx)$ is the neural network model.
Since we hard-constrained the boundary condition, the only loss function will be the residual loss. Hence, there is no weighting hyperparameter.

\begin{table}[htbp]
\centering
\begin{tabular}{|c|c|c|c|c|c|c|}
\hline
 & Metric/Dimension & $10^1$ & $10^2$ & $10^3$ & $10^4$ & $10^5$ \\ \hline\hline
\multirow{5}{*}{Vanilla MC-tfPINN} & Speed & 864it/s & 342it/s & 108it/s & 18it/s & 2it/s \\ \cline{2-7} 
 & Error\_TwoBody $\alpha=0.5$ & 1.86E-3 & 7.27E-3 & 7.41E-3 & 3.56E-3 & 2.17E-3 \\ \cline{2-7} 
 & Error\_TwoBody $\alpha=1.5$ & 1.79E-3 & 7.85E-3 & 8.00E-3 & 4.57E-3 & 2.36E-3 \\ \cline{2-7} 
 & Error\_ThreeBody $\alpha=0.5$ & 4.01E-3 & 4.94E-4 & 4.32E-4 & 4.54E-4 & 5.69E-4 \\ \cline{2-7} 
 & Error\_ThreeBody $\alpha=1.5$ & 3.97E-3 & 5.02E-4 & 4.44E-4 & 5.02E-4 & 5.58E-4 \\ \hline\hline
\multirow{5}{*}{Improved MC-tfPINN} & Speed & 1036it/s & 468it/s & 117it/s & 21it/s & 3it/s \\ \cline{2-7} 
 & Error\_TwoBody $\alpha=0.5$ & 2.03E-3 & 8.92E-3 & 6.57E-3 & 3.64E-3 & 2.09E-3 \\ \cline{2-7} 
 & Error\_TwoBody $\alpha=1.5$ & 1.97E-3 & 8.08E-3 & 6.42E-3 & 4.08E-3 & 2.11E-3 \\ \cline{2-7} 
 & Error\_ThreeBody $\alpha=0.5$ & 4.02E-3 & 5.62E-4 & 4.00E-4 & 4.81E-4 & 5.09E-4 \\ \cline{2-7} 
 & Error\_ThreeBody $\alpha=1.5$ & 3.81E-3 & 5.59E-4 & 4.25E-4 & 4.97E-4 & 5.28E-4 \\ \hline
\end{tabular}
\caption{Results for Monte Carlo tempered fractional PINN on various inseparable exact solutions with a tempered fractional Poisson equation. The error is stable, around 1e-3 to 1e-4, and the cost is sublinear. The two models perform similarly so we do not bold the better result.}
\label{tab:tfpinn_tfpoisson_results}
\end{table}

The results are shown in Table \ref{tab:tfpinn_tfpoisson_results} with different methods' speed (iteration per second) and relative $L_2$ error. This test case is relatively simple. Thus, both variants of tempered fractional PINN yield low errors ranging from 1e-3 to 1e-4. Additionally, we validated that the improved version exhibits faster convergence. In this scenario, the forcing term of the tempered fractional PDE is not analytical, hence we resort to Monte Carlo to unbiasedly estimate the forcing term. In other words, the forcing term used during the optimization of MC-tfPINN is noisy. Despite the presence of noise in the forcing term, both variants of tfPINN still yield low errors, demonstrating the robustness of our approach.

\subsection{Time-Dependent Tempered Fractional Diffusion Equation}
We demonstrate how our MC-tfPINN can deal with a time-dependent tempered fractional diffusion equation:
\begin{align}
\frac{\partial^{\gamma} u(\bx, t)}{\partial t^\gamma} + c \left(-\Delta_\lambda\right)^{\alpha/2}u(\bx, t) + \bv \cdot \nabla_{\bx}u(\bx, t) &= f(\bx, t), \quad (\bx, t) \in \Omega \times (0, T].\\
u(\bx, t=0)&= 0, \quad \bx \in \Omega.\\
u(\bx, t) &= 0, \quad (\bx, t) \in \Omega^c \times (0, T].
\end{align}
The fabricated solution is highly coupled:
\begin{equation}\label{eq:sol_2body_t}
u_{\text{two-body}}(\bx, t) = \text{ReLU}\left(1 - \Vert \bx \Vert_2^2\right)\left(\sum_{i=1}^{d-1}  c_i \sin\left[t \cdot (\bx_i +\cos(\bx_{i+1})+\bx_{i+1}\cos(\bx_i))\right]\right),
\end{equation}
where $c_i$ is randomly drawn from a unit Gaussian.
The zero boundary and initial conditions pose additional challenges in residual fitting on the given complicated forcing term $f(\bx, t)$ since they will not leak any information about the exact solution. The terminal time is set to $T=1$, the time-fractional order is set to $\gamma=0.5$, and the velocity $\bv \in \mathbb{R}^d$ is uniformly drawn from $[0, 1]$. 
For other coefficients in the PDE operator, we choose the following combinations: (1) $\alpha = 0.5, \lambda=1$; (2) $\alpha = 1.5, \lambda=1$; (3) $\alpha = 0.5, \lambda=2$; (4) $\alpha = 0.5, \lambda=0.5$. We want to test if our MC-tfPINN and its improved version can deal with diverse settings with different PDE operator coefficients.

Here are the implementation details. The model is a 4-layer fully connected network with 128 hidden units activated by $\tanh$, which is trained via Adam \cite{kingma2014adam} for 100K epochs, with an initial learning rate 1e-3, which linearly decays to zero at the end of the optimization. We select 100 random residual points at each Adam epoch and 20K fixed testing points uniformly from the unit ball. We adopt the following model structure to satisfy the zero boundary condition with hard constraint and to avoid the boundary and initial losses \cite{lu2021physics}
$
t \text{ReLU} (1 - \Vert\bx\Vert_2^2) u_\theta(\bx,t),
$
where $u_\theta(\bx,t)$ is the neural network model.
Since we hard-constrained the boundary and initial conditions, the only loss function will be the residual loss. Hence, there is no weighting hyperparameter.

\begin{table}[htbp]
\centering
\begin{tabular}{|c|c|c|c|c|c|c|c|}
\hline
Method & $\alpha$ & $\lambda$ & Dimension & 10 & 100 & 1,000 & 10,000 \\ \hline
\multirow{2}{*}{Vanilla MC-tfPINN} & \multirow{2}{*}{0.5} & \multirow{2}{*}{1} & Speed & 182it/s & 161it/s & 20it/s & 3it/s \\ \cline{4-8} 
 &  &  & Rel. $L_2$ & 2.49E-02 & 4.94E-02 & 6.91E-02 & 7.13E-02 \\ \hline
\multirow{2}{*}{Improved MC-tfPINN} & \multirow{2}{*}{0.5} & \multirow{2}{*}{1} & Speed & \textbf{202it/s} & \textbf{189it/s} & \textbf{21it/s} & \textbf{4it/s} \\ \cline{4-8} 
 &  &  & Rel. $L_2$ & \textbf{1.63E-02} & \textbf{2.69E-02} & \textbf{3.73E-02} & \textbf{4.43E-02} \\ \hline
\hline
Method & $\alpha$ & $\lambda$ & Dimension & 10 & 100 & 1,000 & 10,000 \\ \hline
\multirow{2}{*}{Vanilla MC-tfPINN} & \multirow{2}{*}{1.5} & \multirow{2}{*}{1} & Speed & 182it/s & 161it/s & 20it/s & 3it/s \\ \cline{4-8} 
 &  &  & Rel. $L_2$ & \multicolumn{1}{l|}{2.97E-02} & \multicolumn{1}{l|}{3.33E-02} & \multicolumn{1}{l|}{4.79E-02} & \multicolumn{1}{l|}{4.94E-02} \\ \hline
\multirow{2}{*}{Improved MC-tfPINN} & \multirow{2}{*}{1.5} & \multirow{2}{*}{1} & Speed & \textbf{202it/s} & \textbf{189it/s} & \textbf{21it/s} & \textbf{4it/s} \\ \cline{4-8} 
 &  &  & Rel. $L_2$ & \textbf{2.10E-02} & \textbf{2.26E-02} & \textbf{4.79E-02} & \textbf{4.61E-02} \\ \hline
 \hline
Method & $\alpha$ & $\lambda$ & Dimension & 10 & 100 & 1,000 & 10,000 \\ \hline
\multirow{2}{*}{Vanilla MC-tfPINN} & \multirow{2}{*}{0.5} & \multirow{2}{*}{2} & Speed & 182it/s & 161it/s & 20it/s & 3it/s \\ \cline{4-8} 
 &  &  & Rel. $L_2$ & \multicolumn{1}{l|}{4.10E-03} & \multicolumn{1}{l|}{3.67E-02} & \multicolumn{1}{l|}{5.39E-02} & \multicolumn{1}{l|}{5.53E-02} \\ \hline
\multirow{2}{*}{Improved MC-tfPINN} & \multirow{2}{*}{0.5} & \multirow{2}{*}{2} & Speed & \textbf{202it/s} & \textbf{189it/s} & \textbf{21it/s} & \textbf{4it/s} \\ \cline{4-8} 
 &  &  & Rel. $L_2$ & \textbf{4.04E-03} & \textbf{4.89E-03} & \textbf{4.69E-02} & \textbf{4.59E-02} \\ \hline
 \hline
Method & $\alpha$ & $\lambda$ & Dimension & 10 & 100 & 1,000 & 10,000 \\ \hline
\multirow{2}{*}{Vanilla MC-tfPINN} & \multirow{2}{*}{0.5} & \multirow{2}{*}{0.5} & Speed & 182it/s & 161it/s & 20it/s & 3it/s \\ \cline{4-8} 
 &  &  & Rel. $L_2$ & \multicolumn{1}{l|}{9.31E-02} & \multicolumn{1}{l|}{8.98E-02} & \multicolumn{1}{l|}{6.91E-02} & \multicolumn{1}{l|}{7.70E-02} \\ \hline
\multirow{2}{*}{Improved MC-tfPINN} & \multirow{2}{*}{0.5} & \multirow{2}{*}{0.5} & Speed & \textbf{202it/s} & \textbf{189it/s} & \textbf{21it/s} & \textbf{4it/s} \\ \cline{4-8} 
 &  &  & Rel. $L_2$ & \textbf{1.61E-03} & \textbf{1.07E-02} & \textbf{4.56E-02} & \textbf{4.59E-02} \\ \hline
\end{tabular}
\caption{Results of MC-tfPINN on a time-dependent tempered fractional diffusion equation. Improved MC-tfPINN is more efficient and has lower error thanks to its accurate quadrature. Faster speed and lower error are bold.}
\label{tab:tfpinn_time}
\end{table}

The results are presented in Table \ref{tab:tfpinn_time}. Firstly, this test case is notably more challenging compared to the previous one, primarily due to the time-fractional derivative.
Additionally, the exact solution in this case is highly coupled in both time and spatial $\bx$, rendering it exceedingly complex and inseparable.
However, our two variants of tfPINN still manage to yield around 5\% error across all dimensions, with speed growing sublinearly as the dimension increases and the error increasing slowly with the dimension. The speeds are the same if the PDE dimensionality and the model are the same since the PDE operator coefficient does not affect the speed. Furthermore, we reaffirm that improved MC-tfPINN with quadrature always trains faster and achieves smaller errors. In conclusion, we extensively tested MC-tfPINN on tempered fractional PDEs with different PDE operator parameters, demonstrating that tfPINN is fairly stable with respect to the parameters alpha (fractional order) and lambda (tempering factor) of the PDE operator. This showcases the stability and scalability of our proposed MC-tfPINN across different dimensions and various PDE operator parameter settings.

\subsection{Inverse Tempered Fractional PDEs}
We will test the two-body interaction exact solution in equation (\ref{eq:sol_2body}) with the tempered fractional Laplacian equation $\left(-\Delta_\lambda\right)^{\alpha / 2} u = f$ in the unit ball.
Unlike the forward problem, we know the values of the exact solution on several points and wish to infer the fractional order $\alpha$ or the tempering factor $\lambda$. The implementation details on optimizing over $\alpha$ or $\lambda$ are presented in Appendix \ref{sec:mctfpinn_inverse_detail}. Specifically, we consider the following two inverse problems.

\textbf{Problem 1: Inverse problem of $\alpha$}. We fix $\lambda = 1$ and initial the unknown fractional order as $\alpha = 1.5$. The exact fractional order $\alpha$ is either 1.4 or 1.6.

\textbf{Problem 2: Inverse problem of $\lambda$}. We fix $\alpha = 1.5$ and initial the unknown tempering factor as $\lambda = 1$. The exact $\lambda$ is set to 2.

The forcing terms for training residual points are computed with a Monte Carlo simulation with 1024 samples, while the Monte Carlo/quadrature sample size for the neural network is 128 during training. $\epsilon = 1e-6$ for MC-tfPINN to truncate $r$ to prevent numerical instability.

Here are more implementation details. The model is a 4-layer fully connected network with 128 hidden units activated by $\tanh$, which is trained via Adam \cite{kingma2014adam} for 100K epochs, with an initial learning rate 1e-3, which linearly decays to zero at the end of the optimization. We select 100 random residual points and data points at each Adam epoch. We adopt the following model structure to satisfy the zero boundary condition with hard constraint and to avoid the boundary loss \cite{lu2021physics}
$
\text{ReLU} (1 - \Vert\bx\Vert_2^2) u_\theta(\bx),
$
where $u_\theta(\bx)$ is a fully connected neural network.
Since we hard-constrained the boundary condition, the loss functions will be the residual loss and the data loss in the inverse problem and the weights are one.

\begin{table}[]
\centering
\begin{tabular}{|c|c|c|c|c|c|c|}
\hline
Method & Problem & Dimension & Exact $\alpha$ & Identified $\alpha$ & Rel $L_1$ & Time \\ \hline \hline
MC-tfPINN & 1 & 100 & 0.4 & 0.3902741 & 2.43\% & 6min \\ \hline
Improved MC-tfPINN & 1 & 100 & 0.4 & \textbf{0.398084} & \textbf{0.48\%} & \textbf{5min} \\ \hline \hline
MC-tfPINN & 1 & 1000 & 0.4 & 0.3905588 & 2.36\% & 17min \\ \hline
Improved MC-tfPINN & 1 & 1000 & 0.4 & \textbf{0.3989002} & \textbf{0.27}\% & \textbf{16min} \\ \hline\hline
MC-tfPINN & 1 & 10000 & 0.4 & 0.3907224 & 2.32\% & 103min \\ \hline
Improved MC-tfPINN & 1 & 10000 & 0.4 & \textbf{0.3988998} & \textbf{0.28\%} & \textbf{102min} \\ \hline\hline
MC-tfPINN & 1 & 100 & 0.6 & 0.5896347 & 1.73\% & 6min \\ \hline
Improved MC-tfPINN & 1 & 100 & 0.6 & \textbf{0.5983624} & \textbf{0.27\%} & \textbf{5min} \\ \hline\hline
MC-tfPINN & 1 & 1000 & 0.6 & 0.5901451 & 1.64\% & 17min \\ \hline
Improved MC-tfPINN & 1 & 1000 & 0.6 & \textbf{0.5992336} & \textbf{0.13\%} & \textbf{16min} \\ \hline\hline
MC-tfPINN & 1 & 10000 & 0.6 & 0.5902596 & 1.62\% & 103min \\ \hline
MC-tfPINN & 1 & 10000 & 0.6 & \textbf{0.5991701} & \textbf{0.14\%} & \textbf{102min} \\ \hline\hline
Method & Problem & Dimension & Exact $\lambda$ & Identified $\lambda$ & Rel $L_1$ & Time \\ \hline\hline
MC-tfPINN & 2 & 100 & 2 & 2.005856 & 0.29\% & 5min \\ \hline
Improved MC-fPINN & 2 & 100 & 2 & \textbf{2.00189} & \textbf{0.09\%} & \textbf{4min} \\ \hline\hline
MC-tfPINN & 2 & 1000 & 2 & 2.004311 & 0.22\% & 18min \\ \hline
Improved MC-tfPINN & 2 & 1000 & 2 & \textbf{2.000428} & \textbf{0.02\%} & \textbf{17min} \\ \hline\hline
MC-tfPINN & 2 & 10000 & 2 & 2.00E+00 & 0.20\% & 116min \\ \hline
Improved MC-tfPINN & 2 & 10000 & 2 & \textbf{2.000902} & \textbf{0.05\%} & \textbf{114min} \\ \hline
\end{tabular}
\caption{Results for Monte Carlo tempered fractional PINN on inverse problems, where we report the problem setting and the relative $L_1$ error for the identified PDE operator coefficient. We also present the total running time for the algorithms. We bold the best results.}
\label{tab:inverse_tfpinn_tfpoisson_results}
\end{table}

Table \ref{tab:inverse_tfpinn_tfpoisson_results} presents the identified PDE operator coefficient value and its corresponding relative $L_1$ error given the exact coefficient value.  MC-tfPINN can satisfactorily solve both forward and inverse PDE problems, demonstrating its versatility. In addition, we further prove that the improved MC-tfPINN using quadrature can be trained faster and yield better results thanks to its improved accuracy and efficiency.

\section{Summary}
We extend the application of physics-informed neural networks (PINNs) to high-dimensional fractional and tempered fractional partial differential equations (PDEs) by developing the Monte Carlo tempered fractional PINN (MC-tfPINN) based on the previous work Monte Carlo fPINN \cite{guo2022monte} (MC-fPINN).
Furthermore, this research enhances MC-fPINN and MC-tfPINN with Gaussian quadrature to replace traditional Monte Carlo methods when computing the 1D radial singular integral within their PDE operators. This new approach significantly reduces computational errors and increases convergence speed by addressing the limitations associated with high variance and hyperparameter sensitivity in previous models MC-fPINN \cite{guo2022monte}, which rely on a pure Monte Carlo simulation. 
The proposed methods are validated across a range of forward and inverse problems of (tempered) fractional PDEs, demonstrating their superior accuracy and efficiency in handling up to 100,000 dimensions, a marked improvement over existing methods.
We also show that the improved methods based on quadrature can consistently outperform the vanilla models based on Monte Carlo solely in terms of speed and accuracy.

\newpage

\appendix
\section{Unit Test on Quasi-Monte Carlo (QMC)}\label{exp:mc_tfpinn_qmc}
We conduct a unit test on the accuracy of different fPINN simulations in estimating the fractional Laplacian. Specifically, we consider the following three models.
\begin{enumerate}
\item The vanilla MC-fPINN \cite{guo2022monte}. It uses Monte Carlo for the integral of both $r$ and $\boldsymbol{\xi}$.
\item Our improved MC-fPINN. It uses Gauss-Jacobi quadrature for the 1D integral of $r$ and conventional Monte Carlo for the integral of $\boldsymbol{\xi}$.
\item QMC-based improved fPINN. It uses Gauss-Jacobi quadrature for the 1D integral of $r$ and Quasi-Monte Carlo for the integral of $\boldsymbol{\xi}$.
\end{enumerate}
We will test the following exact solution with the fractional Laplacian equation $(-\Delta)^{\alpha / 2} u = f$ in the unit ball:
\begin{equation}
u_{\text{exact}}(\bx) = \left(1 - \Vert \bx \Vert^2\right)^{\alpha/2}\left(c_{1, 0} + \sum_{i=1}^d c_{1,i}\bx_i\right) + \left(1 - \Vert \bx \Vert^2\right)^{1 + \alpha/2}\left(c_{2, 0} + \sum_{i=1}^d c_{2,i}\bx_i\right),
\end{equation}
where $c_{1, i}$ and $c_{2, i}$ are independently sampled from a unit Gaussian. The exact forcing term $f$ can be computed analytically based on Dyda \cite{Dyda+2012+536+555}.
We conducted the 100D and 1000D test cases and repeated the experiment ten times using ten different random seeds.
We report the relative $L_2$ error for fractional Laplacian forcing term estimation in Table \ref
{tab:mctfpinn_unittest}. 
We randomly sampled 20K test points within the unit ball to test different methods' accuracy in estimating the forcing term. We choose 64 points for the Monte Carlo/quadrature batch size. $\epsilon$ for preventing numerical issues and truncating the sampled $r$ from Beta distribution in MC-fPINN is chosen as 1e-6. $\alpha$ is set to 1.5.
First, our improved MC-fPINN outperforms MC-fPINN \cite{guo2022monte} in lower error and variance thanks to the accurate quadrature. Besides, improved MC-fPINN and QMC-fPINN perform similarly, i.e., the QMC in the latter does not necessarily reduce variance and error in high dimensions. These results validate the effectiveness of the key quadrature component in our improved MC-fPINN. 
We observe from the table that QMC in the integral of $\boldsymbol{\xi}$ does not improve the accuracy. 

\begin{table}[htbp]
\centering
\begin{tabular}{|c|c|c|}
\hline
Method/Dimension & 100D & 1,000D \\ \hline
MC-fPINN & 1.32E-2$\pm$1.23E-2 & 4.05E-3$\pm$1.035E-3 \\ \hline
Improved MC-fPINN & \textbf{1.03E-2$\pm$6.89E-3} & \textbf{4.02E-5$\pm$3.83E-5} \\ \hline
QMC-fPINN & 1.24E-2$\pm$7.02E-3 & 6.21E-5$\pm$1.66E-5 \\ \hline
\end{tabular}
\caption{Results for the unit test on QMC.}
\label{tab:mctfpinn_unittest}
\end{table}

\section{Implementation Detail on Inverse Problem}\label{sec:mctfpinn_inverse_detail}
Application to inverse problems requires careful handling of the fractional order. Unlike the forward problem, we need to optimize the PDE coefficient in the inverse problem.

\subsection{Time-fractional Derivative}
\label{sec:time-fractional-represenation}
By the representation of the Caputo time-fractional derivative \eqref{eq:time-fractional-derivative-equiv} and \eqref{eq:time-fractional-derivative-as-expecation},  we have 
\begin{align*}
\frac{\partial^\gamma f(t)}{\partial t^\gamma} 
&=\gamma t^{1-\gamma}\int^1_0 (\tau)^{-\gamma} \frac{f(t) - f(t - t\tau)}{t\tau}d\tau + \frac{f(t) - f(0)}{t^\gamma}\\
&= \frac{\gamma}{1 - \gamma}t^{1-\gamma}\mathbb{E}_{\tau \sim \text{Beta}(1-\gamma, 1)}\left[\frac{f(t) - f(t - t\tau)}{t\tau}\right]+ \frac{f(t) - f(0)}{t^\gamma}.
\end{align*}
We want to infer $\gamma$ while the Monte Carlo simulation is based on a distribution with $\gamma$ as a parameter. Thus, resampling is required once we have an update on $\gamma$. 
Hence, we assume the unknown order $\gamma \in (\gamma_L, \gamma_H)$, then
\begin{align*}
\frac{\partial^\gamma f(t)}{\partial t^\gamma}
&=\gamma t^{1-\gamma}\int^1_0 (\tau)^{-\gamma} \frac{f(t) - f(t - t\tau)}{t\tau}d\tau + \frac{f(t) - f(0)}{t^\gamma}\\
&=\gamma t^{1-\gamma}\int^1_0 (\tau)^{-\gamma_H} (\tau)^{\gamma_H-\gamma} \frac{f(t) - f(t - t\tau)}{t\tau}d\tau + \frac{f(t) - f(0)}{t^\gamma}\\
&= \frac{\gamma}{1 - \gamma_H}t^{1-\gamma}\mathbb{E}_{\tau \sim \text{Beta}(1-\gamma_H, 1)}\left[(\tau)^{\gamma_H-\gamma}\frac{f(t) - f(t - t\tau)}{t\tau}\right]+ \frac{f(t) - f(0)}{t^\gamma}.
\end{align*}
Then we may optimize over $\gamma$ without resampling since $\gamma_H$ is fixed. We  parameterize the trainable $\gamma$ as follows:
\begin{align}
\gamma = \operatorname{Sigmoid}(\Gamma) \cdot (\gamma_H - \gamma_L) + \gamma_L,
\end{align}
where $\operatorname{Sigmoid}$ controls its output value in $(0, 1)$ and $\Gamma$ is a unconstraint trainable parameter. We initial $\Gamma = 0$ to initialize $\gamma = (\gamma_H + \gamma_L) / 2$ as 
$\operatorname{Sigmoid}(0) = 0.5$.

\subsection{Tempered Fractional Laplacian}
For the tempered fractional Laplacian, we don't need to assume prior knowledge on $\lambda>0$ but require assuming $\alpha \in (\alpha_L, \alpha_H)$. 
By the representation \eqref{eq:fractional-laplacian-as-expectation}, we have
\begin{align*}
&\quad\int_{\mathbb{R}^{d}} \frac{u(\bx) - u(\by)}{\exp\left(\lambda\Vert\bx - \by\Vert\right)\Vert \bx - \by \Vert^{d + \alpha}}d\by \\
&= \frac{1}{2} \int_{\mathcal{S}^{d-1}}\int_{0}^\infty \frac{2u(\bx) - u(\bx - \boldsymbol{\xi} r) - u(\bx + \boldsymbol{\xi} r)}{r^{2}}\exp\left(-\lambda r\right)r^{1-\alpha}drd\boldsymbol{\xi}\\
&= \frac{1}{2} \lambda^\alpha \int_{\mathcal{S}^{d-1}}\int_{0}^\infty \frac{2u(\bx) - u(\bx - \boldsymbol{\xi} r / \lambda) - u(\bx + \boldsymbol{\xi} r / \lambda)}{r^{2} }\exp\left(-r\right)r^{1-\alpha}dr d\boldsymbol{\xi}  \\
&= \frac{1}{2} \lambda^\alpha \int_{\mathcal{S}^{d-1}}\int_{0}^\infty \frac{2u(\bx) - u(\bx - \boldsymbol{\xi} r / \lambda) - u(\bx + \boldsymbol{\xi} r / \lambda)}{r^{2} }\exp\left(-r\right)r^{1-\alpha_H}r^{\alpha_H - \alpha}dr d\boldsymbol{\xi}  \\
&= \frac{1}{2}\lambda^\alpha|\mathcal{S}^{d-1}|\Gamma(2 - \alpha_H) \mathbb{E}_{r, \boldsymbol{\xi}} \left[\frac{2u(\bx) - u(\bx - \boldsymbol{\xi} r / \lambda) - u(\bx + \boldsymbol{\xi} r / \lambda)}{r^{2}}r^{\alpha_H - \alpha}\right],
\end{align*}
where $r$ obeys the Gamma distribution $\operatorname{Gamma}(2 - \alpha_H, 1)$ and $\boldsymbol{\xi}$ is uniformly distributed the $d$-dimensional unit sphere. Hence, we can optimize over $\alpha$ or $\lambda$ without resampling.


\bibliographystyle{plain}
\bibliography{ref}
\end{document}